

\documentstyle{amsppt}
\magnification=1100
\NoBlackBoxes
\def\C{\Bbb C}

\def\P{\Bbb P}

\def\Z{\Bbb Z}
\def\O{\Cal O}
\def\X{\Cal X}

\def\ls{\vskip.25in}
\def\ss{\vskip.15in}

\def\[{\big[}
\def\]{\big]}
\def\V{\bar{V}}

{\tt{000315}}\ss

\centerline{\bf THE DEGREE OF THE DIVISOR OF JUMPING RATIONAL CURVES}
\vskip.25in
\centerline{Z. Ran}
\centerline{Math Dept. UCR, Riverside CA 92521 USA}
\centerline{ {ziv\@ math.ucr.edu}}\ls
\centerline{\bf Abstract}\ss

\font\small=cmr8
{\small {For a semistable reflexive sheaf $E$ of rank $r$ and
$c_1=a$ on $\P^n$ and an integer $d$ such that $r|ad$,
we give sufficient conditions so that the restriction of
$E$ on a generic rational curve of degree $d$ is balanced,
i.e. a twist of the trivial bundle (for instance, if $E$ has balanced
restriction on a generic line, or $r=2$ or $E$ is an exterior
power of the tangent bundle). Assuming this,
we give a formula for the 'virtual degree', interpreted
enumeratively, of the (codimension-1) locus of rational curves of degree
$d$ on which the restriction of $E$ is not balanced,
generalizing a classical formula due to Barth for the
degree of the divisor of jumping lines of a semistable
rank-2 bundle. This amounts to computing a certain determinant line bundle
associated to $E$ on a parameter space for rational curves,
and is closely related to the 'quantum K-theory' of projective
space.}\par

}
\ls
\centerline{\bf{INTRODUCTION}}\ls
Let $E$ be a  reflexive sheaf (or vector bundle)
of rank $r$
on a projective space $\P^n, n\geq 2.$
By a theorem of Grothendieck, the pullback
$E_C$ of $E$ to the nonsingular model of any
rational curve in $\P^n$ can be decomposed as
$$E_C=\bigoplus\limits_{i=1}^{r} \O(k_i),\ \ \ k_1\geq\ldots\geq k_r,$$
where the sequence $(k_1,\ldots ,k_r)$ is uniquely determined
and called the {\it {splitting type}} of $E$ on $C$. By
semi-continuity, splitting type determines a natural locally
closed stratification of any parameter space for rational
curves, which is an important aspect of the geometry of $E$.\par
For $E$ semistable of rank $r=2$,
a well-known theorem of Grauert-M\"ulich
[OSS] says that if $r=2$, the splitting type of
$E$ on a generic line
$L\subset\P^n$ is either $(k,k)\ (c_1(E)$ even) or
$(k,k-1)\ (c_1(E)$ odd).
One of the most important geometric objects associated
to $E$ of rank 2 is the locus of {\it{jumping lines}}, i.e.
the locus of lines $L$ such that $E_L$ is not 
of the above generic form.
This locus is especially important in case $c_1(E)$
is even,
when it is actually a divisor. A well-known formula due
originally to Barth  computes the degree of the divisor
of jumping lines (with a suitable scheme structure).\par
The purpose of this paper is to generalize Barth's formula
to the case of bundles of arbitrary rank $r$ and 
rational curves of any degree $d$ such that
$$r|dc_1(E),\tag 0.1$$ 
essentially under the condition that the
restriction of $E$ on a generic rational curve $C$ of degree
$d$ is balanced, i.e. has splitting type $(k^r),$
so that the locus of rational curves
$C$ for which $E_C$ is unbalanced is of pure codimension
1 and may be called the  {\it{divisor of jumping rational
curves}} of degree $d$ of $E$.\par
Now the balancedness  condition is not satisfied
for all semistable sheaves $E$ even for lines.
However, we will introduce below
a condition (condition AB), 
which states that the restriction of $E$ on a generic line $L$
is 'almost balanced', i.e. has splitting type
 $(k^s,(k-1)^{r-s})$ for some $k,s$,
(so if $r|c_1(E)$ then $E_L$ is actually balanced).
Assuming this, we will then introduce 
a certain transversality condition T and show that with this condition
the restriction $E_C$ to a generic rational curve 
of any degree $d$ is almost balanced, so
that under (0.1) a divisor of jumping rational curves of
degree $d$
may be defined. We will also show that conditions AB and T are
satisfied whenever either $r|c_1(E)$ and AB holds, or $r=2$ or $E$ is
an exterior power of the tangent bundle.\par
Our main result (Theorem 3.1), which assumes
conditions AB and T, computes the 'degree' of the divisor of
jumping rational curves, which is interpreted as the weighted number
of these curves incident to a generic collection of linear
spaces (assuming of course that the total codimension  of the incidence
conditions equals the dimension of the divisor). More precisely,
our formula expresses this degree in terms of some other enumerative
invariants which have been computed before, 
e.g. in [P], [R1],[R2],[R3].\par
Note that even a homogeneous bundle like the tangent
bundle in general admits jumping rational curves of degree $d>n$.
In fact, the divisors of jumping rational curves associated to homogeneous
vector bundles are an interesting class of projectively invariant
divisors on any parameter space of rational curves (be it
the Chow variety considered here, or Kontsevich's space
or whatever).\par
This note was motivated by a talk by Givental [G] on his 
'quantum K-theory', still under construction. This theory
seeks to compute expressions of the form
$$\chi (M_{0,m}(\P^n,d), ev^*(E_1\boxtimes ...\boxtimes E_m)),$$
where $ M_{0,m}(\P^n,d)$ is Kontsveich's space of stable $m$-pointed
rational curves of degree $d$ and $ev: M_{0,m}(\P^n,d)\to (\P^n)^m$
is the evaluation map. As will emerge below, our formula
essentially amounts to a computation of 
$$\chi(M_{0,1}(\P^n,d)_{(A.)}, ev^*(E)),$$ where
$M_{0,1}(\P^n,d)_{(A.)}$ is the normalization of a 1-dimensional subvariety of
$M_{0,1}(\P^n,d)$ defined by incidence to a generic collection
$(A.)$ of linear spaces. By the Riemann-Roch
formula, the latter may also be identified as an intersection number
in $M_{0,m}(\P^n,d)$; for instance in case all the $A_i$ are points,
hence have trivial normal bundle, it
equals
$$\chi (M_{0,m}(\P^n,d),ev_1^*(E),ev_2^*(h^{n})\cup ...\cup ev_m^*
(h^{n})),$$
where $h$ is the hyperplane class, 
and in general the 'relative' Euler characteristic 
may be defined  by
$$\chi (X,E,b)=\int_X \text {ch}(E)\cup \text {td} (T_X)\cup b\ \ ,$$
and coincides with the Euler characteristic of the restriction of $E$
on any smooth subvariety $B$ with fundamental class $b$
and trivial normal bundle (which may or may not exist);
in our case $B=(ev_2\times ...\times ev_m)^{-1}(pt.)$
has trivial normal bundle.\par
 Though our
formula is apparently new and independent of any quantum methods, 
it might in principle become
accessible by Givental's methods at some point. 
Indeed these methods, unlike ours, might yield results
when the jumping locus has codimension $>1.$\ss
The paper is organised as follows. In Sect.1 we study
the Chow compactification of the family of
rational space curves. It does not seem to be generally known 
that this compactification is well-behaved at least in codimension
1 (and perhaps in codimension 2 as well, as long as multiple components
don't appear), so we have given a self-contained treatment here.
The results come out as expected.
In Sect.2 we study qualitatively the restriction of a bundle
on a projective space to rational curves, focusing on criteria 
to ensure that the restriction on a generic rational curve is
balanced. Our enumerative formula is given in Sect.3. The
proof is based, not surprisingly, on the Riemann-Roch formula.

\ls
\centerline{\bf{{1. RATIONAL CURVES}}}\ls
Here we review some qualitative results about families
of rational curves in $\P^n$. See also [R1][R2][R3] and references
therein for details and proofs.
In what follows we fix $n\geq 2$ and denote by $\bar{V}_d$ or $\bar{V}_{d,n}$ the
closure in the Chow variety of the locus of 
irreducible nodal (if $n=2$) or
nonsingular (if $n>2$) rational
curves of degree $d$ in $\P^n$, with the scheme structure
as closure, i.e. the reduced structure
(recall that the Chow form
of a reduced 1-cycle $Z$ is just
the hypersurface in $G(n-2,\P^n)$ consisting
of all linear spaces meeting $Z$). Thus
$\V_d$ is irreducible reduced of dimension
$(n+1)d+n-3$.  Let $A_1, \ldots, A_k$ be a
generic collection of linear subspaces of respective codimensions $a_1,
\ldots, a_k, 2\leq a_i \leq n$ in $\P^n$.  We denote by
$$
B= B_d = B_d (a_{\cdot}) = B_d (A_{\cdot})
$$
the normalization of the locus (with reduced structure)
$$
\{(C, P_1, \ldots, P_k ) \ : \ C \in \bar{V}_d , P_i  \in C \cap A_i, i =1
, \ldots, k \}
$$
which
is also the normalization of its projection to $\bar{V}_d$, i.e. the
locus of degree-$d$ rational curves (and their specializations) meeting
$A_1, \ldots, A_k$.
 We have
$$
\dim B = (n+1) d + (n-3) - \sum (a_i - 1) .\tag 1.1
$$
When $\dim B=0$ we set
$$N_d(a.)=\deg (B). \tag 1.2$$
When  $n=2$, all $a_i=2$ so the $a$'s may be suppressed.
For $n>2$,
$k$ is called the {\it length} of the condition-vector $(a.)$.
The numbers $N_d(a.)$, first
computed in general by Kontsevich and Manin,
 are computed in [R2],[R3] by an elementary method
based on recursion on $d$ and $k$.

Now suppose $\dim B=1$ and let
$$
\pi : X \to B \tag 1.3
$$
be the normalization of the tautological family of rational
curves, and $f: X \to \P^n$ the natural map.
The following summarizes results from [R2][R3] (proved mostly
in the references therein):
\proclaim{Theorem 1.1}(i) X is smooth .\par
(ii) Each fibre $C$ of $\pi$ is either\par
(a) a $\P^1$ on which $f$ is either an immersion with at most one exception
which maps to a cusp
($n=2$) or an embedding ($n>2$); or\par
(b)  a pair of $\P^1$'s meeting transversely once, on which $f$ is an 
immersion with nodal image ($n=2$) or an embedding ($n>2$); or\par 
(c)  if $n=3$, a $\P^1$ on which $f$ is a degree-1 immersion such that $f(\P^1)$ has
a unique singular point which is an ordinary node.\par
(iii) If $n>2$ then $\bar{V}_{d,n}$ is smooth along the image $\bar{B}$ of
$B$, and $\bar{B}$ is smooth except, in case some $a_i=2,$ for
ordinary nodes corresponding to curves meeting some
$A_i$ of codimension 2 twice. 
If $n=2$ then $\bar{V}_{d,n}$ is smooth
in codimension 1 except for a cusp along the cuspidal locus and
normal crossings along the reducible locus, and $\bar{B}$ has the
singularities induced from $\bar{V}_{d,n}$ plus ordinary nodes
corresponding to curves with a node at some $A_i$, and no other
singularities.
\endproclaim
The basic idea is that one can first get a handle on what curves
occur in $\bar{B}$ by the standard technique of semistable
reduction (actually, normalization
is sufficient) plus dimension counting (as, e.g. in Harris' work
on the Severi problem); then doing the deformation theory for the curves
which do occur is easy enough. For the convenience of the reader
we will give a complete proof for $n>2.$\ls
Let $H_d$ denote the (scheme-theoretic) closure in the Hilbert
of the family of nonsingular rational curves of degree $d$, and
$H_d^0\subset H_d$ be the open subset of reduced curves with
normal crossings, which is well known to be smooth (see also below).
Let $$\pi_d:X_d\to \bar{V}_d,\ \ \ \pi_{H_d}:X_{H_d}\to H_d$$ 
be the universal cycle (resp. universal curve). There is a 
natural morphism (cf. [K])
$$c:H_d\to\bar{V}_d$$
which assigns to a curve $C$ its Chow form, which is the divisor on
the Grassmannian $G(n-2,\P^n)$ consisting of subspaces meeting $C$,
with multiplicities if $C$ is not generically reduced.
Thus $\bar{V_d}$ is naturally embedded in a projective
space parametrizing a suitable linear system on the Grassmannian. Clearly 
$c$ is one to one on the subset $H_d^1\subset H_d$ consisting of reduced curves.
In fact, more is true:
\proclaim{Lemma 1.2} c is unramified at any reduced subscheme.\endproclaim
\demo{proof} Let $C$ be reduced and pick a nonzero $v\in T_{[C]}H_d.$
Let $p\in C$ be a general point where $v_p\neq 0,$ and let $L$ be
a general $(n-2)-$ plane through $p$, which is also a general point in some component
of $c(C).$ Let $\tilde{L}$ be the lift of $L$ to $\C^{n+1},$
and $\tilde(p)\subset\tilde(L)$ the 1-dimensional subspace
lifting $p$. Then the
tangent space to the Grassmannian at $L$ may be identified with
$${\text {Hom}} (\tilde{L}, \C^{n+1}/\tilde{L}),$$
while the tangent space to the divisor $c(C)$ is
$$\{\phi\in {\text {Hom}} (\tilde{L}, \C^{n+1}/\tilde{L}):
 \phi (\tilde{p})\equiv 0\mod T_pC\}.$$
Choosing $L$ general through $p$, we can arrange that $$v_p\not\in 
<T_pC,T_pL>,$$ and it follows that as $C$ moves infinitesimally 
according to $v$, we can move $L$ preserving incidence to $C$ and 
going outside of $c(C)$, so $d_{[C]}(v)\neq 0.$ Note that the Lemma 
and the proof are valid for pure-dimensional subschemes of any 
dimension.\qed\enddemo 
Now the basic codimension-1 dimension 
counting result for $\V_d$ is the following 
\proclaim{Proposition 
1.3} Let $W\subset \V_d$ be any codimension-1 subvariety and 
$[C]\in W$ a general curve. Then $C$ is either\ss (i) a smooth 
embedded $\P^1$; or\ss (ii) a pair of smooth embedded $\P^1$'s 
meeting transversely at one point; or\ss (iii) only if $n=3$, an 
irreducible immersed rational curve with one normal crossing.\ss 
Moreover $\V_d$ is smooth at $[C]$ in each case and has tangent 
space $H^0((I_C/I_C^2)^*)$ in case (i) or (ii), or $H^0(N_f)$, 
where $f:\P^1\to C$ is the normalization, in cases (i) or 
(iii).\endproclaim
 \demo{proof} We will use the following variant 
of Kleiman transversality: 
\proclaim{Lemma 1.4} Let $\{Z_s:s\in 
S\}$ be a family of $k$-cycles in $\P^n$ which is 
$PGL_n$-equivariant (i.e. $S$ is $PGL_n$-invariant and
$gZ_s=Z_{gs}$), 
and let $U\subset\P^n$ be a (purely) 
codimension-$c$ subvariety, $c>k$. Then the locus $S_U:=\{s\in 
S:Z_s\cap U\neq\emptyset\}$ is of codimension $c-k$ in 
$S$.\endproclaim \demo{proof} Take any $s\in S$ and any $(c-k-1)-$ 
dimensional subvariety $Q\subset PGL_n$. Kleiman transversality 
says that $$\overline{(\bigcup\limits_{g\in Q} gZ_s)}\cap g_0U=\emptyset$$ for 
general $g_0\in PGL_n$, hence 
$$\overline{(\bigcup\limits_{g\in Q} 
g_0^{-1}gZ_s)}\cap U=\emptyset.$$ This easily implies that the 
intersection of $S_U$ with the PGL$_n$-orbit of $Z_s$ is of 
codimension $c-k$. Thus $S_U$ meets every orbit in codimension 
$c-k$, and it follows easily that $S_U$ is of codimension 
$c-k$.\qed\enddemo \ss
 We now return to the proof of the Proposition. If 
$W$ fails to be $PGL_n$-invariant, then its general element $[C]$ 
is general in $\V_d$, hence smooth. Hence we may assume $W$ is 
$PGL_n$-invariant. By semistable reduction, there exists a family 
$$Y\to T$$ with general fibre $\P^1$ and special fibre $$Y_0=\bigcup 
Y_{0,i}$$ with normal crossings, and with a surjective map 
$h:Y_0\to C$. Set $$h_i=h|_{Y_{0,i}}, h_{i*}[Y_{0,i}]=m_iC_i\  
{\text{or}}\  0, 
d_i=\deg (C_i), k=\#\{i:m_i>0\}.$$ We may assume $m_1>0$ and that 
$C_1$ is non-disconnecting (i.e.$\bigcup\limits_{i>1}C_i$
is connected). Then from Lemma 1.4 it follows that 
$$(n+1)d+n-4=\dim W=\dim\{C\}\leq 
\dim\{C_1\}+\dim\{\bigcup\limits_{i>1}C_i\} -n+2$$ $$\leq 
(n+1)d_1+n-3+\dim\{\bigcup\limits_{i>1}C_i\} -n+2\leq .\ .\ .$$ $$\leq 
(n+1)\sum d_i+k(n-3)-(k-1)(n-2)=(n+1)\sum d_i+n-3-(k-1).$$ It 
follows at once that $\sum d_i=d$, so that all nonzero $m_i$ are 
equal to 1, and that $k=1$ or 2, and in the latter case all 
inequalities above are equalities.\ss 
Now suppose $k=2$, so 
$$C=C_1\cup_p C_2.$$ Then $(C_1,C_2)$ must be a general point in the 
locus of intersecting curves in $\V_{d_1}\times\V_{d_2}$, which 
has codimension $n-2$, and it follows easily that $C_1$ and $C_2$ 
are smooth; moreover as the stabiliser in $PGL_n$ of a point 
$p\in\P^n$ acts transitively on $\P^n-p$ and on $T_p\P^n$, it 
follows easily that $C_1$ and $C_2$ meet transversely once. Now as 
$C$ is a locally complete intersection, the tangent space to $H_d$ 
at $[C]$ is given by $$H^0(N),\ N=(I_C/I_C^2)^*.$$ On the other hand 
the tangent space to the locally trivial deformations is given by 
$$H^0(N'),$$ where $N'$ is the kernel of a map $$N\to \C_p$$ to a 
skyscraper sheaf at $p$, which assigns to a deformation of $C$ the 
corresponding deformation of the germ $(C,p)$, which is just an 
ordinary node (cf. [S]). Note that $N'$ may also be identified as 
the kernel of the natural map $$N_{C_1}\oplus N_{C_2}\to 
T_p\P^n/(T_pC_1+T_pC_2)\ ,$$
$$(v_1,v_2)\mapsto v_{1,p}-v_{2,p},$$ hence clearly $H^1(N')=0$. It follows 
easily as in [S] that the germ of $H_d$ at $[C]$ is smooth and 
maps unramifiedly to the deformation space of $(C,p)$. 
Since the total space of the 
versal deformation of $(C,p)$ is just given by $xy=t$, it is 
smooth, and it follows that the universal curve $X_d$ is smooth 
along its fibre over $[C]$.\ss 
It remains to consider the case 
where $C$ is irreducible of degree $d$, hence given by projecting 
a rational normal curve $$C_d\subset\P^d$$ from a center $$M=\P^{d-n-1},
\  M\cap C_d=\emptyset.$$
Recall that an arbitrary length-$k$ subscheme of $C_d$ spans a
$\P^{k-1}$, so if $C$ is not embedded then for any
$(C_d,M)$ such that $C= \text{proj}_M(C_d)$,
$M$ must meet some $k$-secant $\P^{k-1}$ to $C_d$ in a $\P^{k-2},
\ k\geq 2.$
By a straightforward dimension count, the family of pairs $(C_d,M)$
with the latter property, for any fixed $k$, is of codimension
$$1+(k-1)(n-2)$$ and, under the  mapping 
$$(C_d,M)\mapsto \text{proj}_M(C_d)$$
maps to a subfamily of the same codimension in $\V_d$.
Since $n>2,$ this codimension could equal 1 only if $n=3,\ k=2,$
 hence the set of $C$'s which are not embedded
has codimension 1 only if $n=3$. We can see similarly that if $n=3$
a general nonembedded $C$ has one normal crossing and corresponds
to $M$ meeting the secant variety of $C_d$ in one general point.\ss
Now for any irreducible $C$ of degree $d\geq 3$, pick general points
$$p_1,p_2,p_3\in C$$ and transverse hyperplanes $H_i\ni p_i$, and
let $$f:\P^1\to C$$ be the normalisation and set $p_i'=f^{-1}(p_i)$. 
Consider a space $D$ parametrising
deformations $f'$ of $f$ so that $$f'(p'_i)\in H_i.$$Then
clearly $$T_fD=H^0(T')$$ where $$T'\subset f^*(T_{\P^n})$$
is the (full-rank) subsheaf of vector fields tangent to $H_i$
at $p'_i,\ i=1,2,3.$ As $$f^*(T_{\P^n})(-p'_1-p'_2-p'_3)\subset T',$$
clearly $H^1(T')=0$, hence $D$ is smooth at $f$. Moreover the natural
map 
$$D\to \V_d$$$$f'\mapsto f'(\P^1)$$
 is clearly one to one and is unramified at $f$ by
an evident variant of Lemma 1.2. Therefore $\V_d$ is smooth at $[C]$
(which clearly implies that the normalization of the total space of the universal
cycle, which  has a smooth fibre over $[C]$, is itself smooth
along this fibre). Also, it is clear that the natural map
$T'\to N_f$ induced an isomorphism $$H^0(T')\to H^0(N_f).$$
This completes the proof of Proposition 1.3.\qed

\enddemo\ss
We can now complete the proof of Theorem 1.1. Let $$X'_d\to\V_d$$
be the normalization of the universal cycle. Then $X'_d$ is 
irreducible of dimension $(n+1)d+n-2$ and
nonsingular in codimension 1 . Likewise the fibred product
$$(X'_d)^k\to\V_d$$ is irreducible of dimension $(n+1)d+n-3+k$
and nonsingular in codimension 1 (its singularities come
from singularities of $\V_d$ and repeated singular points of 
fibres). Moreover the fibres of  $(X'_d)^k$ over $\V_d$
are obviously $k$-dimensional. Now consider the incidence variety
$$I=\{ (C,p_1,...,p_k,A_1,...,A_k):p_i\in A_i, i=1,...k\}$$
$$\subset (X'_d)^k\times G(n-a_1,\P^n)\times ...\times G(n-a_k,\P^n)\}.$$
This is obviously a fibre bundle over $(X'_d)^k$ , hence irreducible of
dimension $$(n+1)d+n-3+k+\sum a_i(n-a_i),$$ i.e. codimension $\sum a_i$
in the product, and nonsingular in codimension 1. It follows that if
$$\sum (a_i-1)=(n+1)d+n-4$$ then a general fibre $B$ of $I$
over $G(n-a_1,\P^n)\times ...\times G(n-a_k,\P^n)$ is smooth and
1-dimensional and the image $\bar{B}$ of $B$ in $\V_d$ can be made
disjoint from any given codimension-2 subvariety, hence the curves in
$\bar{B}$ are as claimed. The smoothness of $X$ as in Theorem 1.1 can
be proved similarly by considering a suitable incidence variety in
$$(X'_d)^{k+1}\times G(n-a_1,\P^n)\times ...\times G(n-a_k,\P^n)$$
(with incidence conditions on the first $k$ points). 
Finally, as $B$ is smooth, clearly singularities of $\bar{B}$ come
from curves meeting some $A_i$ more than once (or nontransversely),
and the remaining assertions about these singularities can be proved
by a dimension-counting argument on the rational normal curve similar
to the one we did above.
This completes 
the proof of Theorem 1.1.\qed

It follows from Theorem 1.1 that
 we can speak about a 'general reducible boundary curve' of $V_d$
as being a 1-nodal curve $$C_1\cup_pC_2$$ which is either embedded 
as such (if $n>2$) or maps to a reducible nodal curve with one 
distinguished separating node (if $n=2$).\ls 
\subheading{Remark 1.5} Another fact which 
follows easily from the above discussion is that for any 
$$b=(C,p_1,...,p_k)\in B,$$ we can identify the tangent space $T_bB$ 
with $$\{ v\in H^0(N):v_{p_i}\equiv 0 \mod T_{p_i}A_i, 
i=1,...,k\}$$ where $N$ is either $(I_C/I_C^2)^*$ in case (a) or 
(b), or $N_f$ in case (a) or (c). This implies, in the notation of 
[R3],Sect 2, that not only is $B=\bigcap B_i$ but $$T_bB=\bigcap\limits_i 
T_bB_i$$ as well,  so that $B$ is the complete transverse intersection of 
the $B_i$ in $B^+,$ a fact which was implicitly used in the 
computation o the genus of $B$ in [R3].\ls

 \centerline{\bf{2. JUMPING 
RATIONAL CURVES}}\ls Here we discuss some qualitative generalities 
about restriction of vector bundles from $\P^n$ to rational 
curves. See [OSS] for details on vector bundles over projective 
spaces. 

A vector bundle $E_C$ of rank $r$ on a rational curve $C$ is said to be {\it
{almost balanced}} if it can be decomposed
$$E_C\simeq s\O(k)\oplus (r-s)\O(k-1).\tag 2.1$$
In this case the subsheaf $s\O(k)\subseteq E_C$ is well-defined
and determines a canonical 'positive' subspace $V_C(p)\subseteq E(p)
:=E\otimes\O_p$ for any $p\in C.$\par
Now let $E$ be a semistable reflexive sheaf of rank $r$ 
on $\P^n$ and chern class $c_1(E)=a\in\Z $. 
We introduce the following
\proclaim{Condition AB}
The restriction $E_L$
of $E$ on a general line is almost balanced.
\endproclaim
By the Grauert-M\"ulich Theorem, condition AB
is satisfied if $E$ is semistable of rank 2. Also,
this condition is obviously satisfied whenever 
$$E=\bigwedge^mT_{\P^n}$$ for any $m$ (though it fails
for $E={\text {Sym}}^mT_{\P^n}, m>1$).
Assuming this condition holds, we try to
get a similar conclusion
for rational curves of higher degree. To this end, consider the following
\proclaim{Transversality condition T} Given a general point $p\in\P^n$,
an arbitrary subspace $W\subset E(p)$, and a general line $L\ni p$,
the positive subspce $V_L(p)$ is transverse to $W$.\endproclaim
\proclaim{Lemma 2.1} Assuming condition AB, condition T is satisfied 
provided either\par
(i) $a\equiv 0\mod r$; or\par
(ii)$r=2$; or\par
(iii) $E=\bigwedge^mT_{\P^n}$ .\endproclaim
\demo{proof} Case (i) is obvious since then $V_L(p)=E(p).$
In case (ii), if $a$ is odd, it is well known [OSS] that
by semistability of $E$ the
1-dimensional subspace $V_L(p)$ varies with $L$, which is sufficient.
In case (iii) $V_L(p)$ is the space of
'multiples' of $T_L(p)$  and
again the result is clear.\enddemo
Thus, conditions AB and T are both satisfied whenever either
$E$ has generic splitting type $(k^r)$ on lines, or has rank 2
or is an exterior power of the tangent bundle.\par
Now we want to study restrictions of $E$ either to general rational
curves or to general reducible curves in $\bar{V}_{d,n}$. To this end
we make another definition. Consider a 1-nodal curve
$$C=C_1\cup_pC_2$$
with each $C_i\simeq \P^1.$
A bundle $E_C$ on $C$ is said to be {\it{almost balanced}} if each
$E_{C_i}$ is almost balanced and the induced positive subspaces
$$V_1,\ V_2 \subseteq E(p)$$
are in general position. It is an easy execise that in this case
we can write
$$E_C\simeq (\oplus L_i)\oplus (\oplus M_j)\tag 2.2$$
where for some $k$, each $L_i$ (resp. $M_j$) is
a line bundle of  total degree 
$k$ (resp $k-1$). Moreover the 'positive' subsheaf $\oplus L_i\subseteq E_C$
is canonically defined and for any general point $q\in C_1$ or $C_2$
the corresponding subspace $V\subseteq E(q)$ can be identified in
an evident sense with either $V_1+V_2$ or $V_1\cap V_2$.
Note that the decomposition (2.2) implies easily that for any small
deformation $(E_{C'},C')$ of $(E_C,C)$ where $C'$ is a smooth $\P^1$,
$E_{C'}$ is almost balanced.
\proclaim{Proposition 2.2} Let $E$ be a reflexive sheaf on $\P^n$
satisfying conditions AB and T, and let $C$ be either a general element
or a general reducible boundary element of $\bar{V}_d.$ Then
$E_C$ is almost balanced.\endproclaim
\demo{proof} We use induction on $d$, the case $d=1$ being exactly condition AB.
Assuming the assertion holds for $d-1$, specialise a general curve $C$ of
degree $d$ to a general reducible
$C_0=C_1\cup_pL$ with $L$ a line. By property T, clearly
$E_{C_0}$ is almost balanced, hence so is $E_C$ for a general $C$.\par
For a general reducible $C=C_1\cup _pC_2$, we know $E_{C_i}$ is
almost balanced and moreover, as each $C_i$ may be viewed as a
specialisation of a polygon, it follows that the positive subspaces
$V_i\subseteq E(p)$ may be assumed transverse, hence $E_C$ is almost balanced.\qed
\enddemo\ls
\subheading{Example 2.3} We consider the case of the tangent
bundle $E=T_{\P^n}$. For any irreducible rational curve
$$C\to \P^n$$ of degree $d\equiv 0\mod n,$ given by a polynomial
vector $$(f_0,...,f_n)\in\oplus H^0(\O_C(d)),$$ we have
an exact sequence
$$0\to\O_C\to\oplus\O_C(d)\to E_C\to 0\ .$$
Dualising and twisting by $d$, we see that sections
of $E_C^*(d+k)$ correspond to syzygies of degree $k$
among the $f_i,$ i.e.
relations of the form
$$\sum g_if_i=0,\ \ \ \ g_i\in H^0(\O_C(k)).$$
It is immediate from this that $C$ is a jumping curve
iff it admits a syzygy of degree $\frac{d}{n}-1$, while
any curve admits a syzygy of degree $\frac{d}{n}.$ It is also
easy to see that $C$ is a jumping curve iff the ideal generated
by $f_0,...,f_n$ in the homogeneous coordinate ring of $C$ 
fails to contain
$H^0(\O_C(d+\frac{d}{n}-1))$, while this ideal
never contains $H^0(\O_C(d+\frac{d}{n}-2)).$
\ls\ls
\centerline{\bf{3. COUNTING JUMPING RATIONAL CURVES}}\ls

We continue with the notations above and assume that $E$ satisfies
conditions AB and T. Then we may define the {\it{jumping locus}}
$$\Cal J_{d,E}\subset \bar{V}_d,$$
as a set, to be the closure of the set of irreducible reduced $C$ such
that $E_C$ is not almost balanced. Now if it happens that
$$-r<ad\leq 0,\tag 3.1$$
then it is easy to endow $\Cal J_{d,E}$ with a global scheme
structure: namely let
$$\Pi:\X\to V_d$$
be the tautological family and
$$\Cal F:\X\to\P^n$$
the natural map, and let $\Cal J_{d,E}$ be the 
{\it{Fitting subscheme }}
$$Fit_1(R^1\Pi_*(\Cal F^*E)),$$
defined by 1st fitting ideal of $R^1\Pi_*(\Cal F^*E) .$
This what is done in [OSS] for $d=1$. For our purposes however,
the hypothesis (3.1) is too restrictive. Without it
one can still endow $\Cal J_{d,E}$, at least in the event it
had codimension 1,  with a scheme structure 'slice
by slice' , as we now proceed to do.

We now assume that
$$r|ad.\tag 3.2$$
In view of Proposition 2.2 this implies that
for a general $C\in V_d,$, $E_C$ is in fact {\it balanced}.  
Let $\pi:X\to B, f:X\to \P^n$ be as in Sect. 1, and let
$$s_i\subset X, i=1,...,k \tag 3.3$$
be the tautological section corresponding to $A_i$.
Let $D=D_{t.}$ be any  divisor of the form
$$D=\sum t_is_i\tag 3.4 $$
such that $\sum t_i=\frac{ad}{r}+1$. Set
$$G=f^*(E)(-D).\tag 3.5$$
Thus for a general fibre $F_b=\pi^{-1}(b)$ we have
$$G_{F_b}=r\O(-1).$$
We define $\Cal J_{d,E,B}$ to be the part of the first Fitting
scheme $Fit_1(R^1\pi_*(G))$ supported in the open subset $B^0\subseteq B$
corresponding to irreducible curves. It is easy to see that
this is independent of the choice of twisting divisor $D$.
In particular, taking $D=(ad+r)s_i,$ our scheme structure
coincides with the natural scheme structure  on  $\Cal J_{d,E,B}$
which defined, at least through codimension 1 
over the locus of curves in $\bar{V}_d$ incident to
$A_i$, by virtue of the existence of a canonical section. 
 Now set
$$J_{d,E}(a.)=c_1(\Cal J_{d,E,B}).\tag 3.6 $$
This evidently depends only on the $(a.)$, and it is this
number that we will compute.\par
To state our formula conveniently we introduce
some objects from [R2][R3]. Set
$$m_i=m_i(a.)=-s_i^2, i=1,...,k.\tag 3.7 $$
Note that if $a_i=a_j$ then $m_i=m_j.$ In
particular for $n=2$ they are all equal.
It is shown in [R2][R3] that these numbers can all be
computed recursively in terms of data of lower degree $d$ and
(for $n>2$) lower length. For instance for $n=2$ we have
$$
2m_1 = \sum_{d_1 + d_2 = d} N_{d_1} N_{d_2} d_1 d_2 {\binom{3d-4}{3d_1-2}}.\tag 3.8
$$
Note that

$$
s_i . s_j = N_d(...,a_i+a_j,...,\hat a_j,...), i \neq j,\tag 3.9
$$
so this number may be considered known. Hence $D^2$ may be
considered known.
Also, let $R_i $ be the sum of all
fibre components not meeting $s_i$ . We have
$$R_i.s_j=m_i+m_j+2s_i.s_j \tag 3.10 $$
so this number is computable as well
(actually we showed in [R2] that this is
computable in terms of lower-degree data,
and the $m_i$ were computed from that).\par
Next, set $$L=f^*(\O(1)),$$ and note that
$$L^2=N_d(2,a.),\  L.s_i=0,\ i=1,...,k.\tag 3.11$$
Also, it is easy to see as in [R3] that
$$L.R_i=\sum_{d_1+d_2=d}\binom{3d-1}{3d_1-1}d_1d_2^2N_{d_1}N_{d_2},\  n=2$$
$$
L.R_i = \sum d_2N_{d_1} (a_{\cdot}^1 , a_i,{n_1} )
N_{d_2} (a^2_{\cdot},
 {n_2} ), \ n>2,
\tag 3.12
$$
the summation for $n>2$ being over all $d_1 + d_2 = d, n_1 + n_2 = n $
and all decompositions
$$A_{\cdot} = (A_i) \coprod (A_{\cdot}^1) \coprod (A_{\cdot}^2)
$$ (as unordered sequences or partitions).\par
Recall that
the geometric genus  $g=g(B)$ was computed in [P][R2]
for $n=2$ and in [R3] in general.
Next, define for any index-set $I\subseteq \{ 1,...,k\}$ of cardinality $|I|$
$$t_I=\sum_{i\in I}t_i,\tag 3.13$$
$$h(a.,t.)=\sum N_{d_1}N_{d_2}\binom{3d-1-|I|}{3d_1-1-|I|}(rt_I-ad_1-r),\  n=2$$
$$h(a.,t.)=\sum N_{d_1}(a_i:i\in I)N_{d_2}(a_i:i\not\in I)(rt_I-ad_1-r),\  n>2\tag 3.14$$
the sum being extended over all $d_1+d_2=d$ and all index-sets $I$
such that $t_I>\lceil\frac{ad_1}{r}\rceil$ (and, for $n=2$,
also $|I|\leq 3d_1-1$).
Now we can finally state our formula.
\proclaim{Theorem 3.1} Notations as above,
assume $E$ satisfies conditions AB and T and set $b=c_2(E)\in\Z$. Then 
the weighted number of jumping rational curves
meeting $A_1,...,A_k$ is given by\ls\noindent
$(3.15)\ \ \ \ \ \ J_{d,E}(a.)=$
$$r(g-1)-\frac{1}{2}((a^2-2b)L^2+rD^2)
+\frac{1}{2}\[(-r)(2g-2-m_1)$$
$$+aR_1.L
+2rs_1.D
-\sum_{i=2}^kt_iR_1.s_i)\]-h(t.,a.)\ \ \ \ \ \ \ \ \hskip1in $$\ls
\endproclaim
\demo{proof} We apply the Riemann-Roch formula in Grothendieck's
form ([F],15.2) (though the Hirzebruch form would have worked too)
to the vector bundle $G$ and the map $\pi :X\to B$. Clearly
$\pi_*(G)=0$ while $R^1\pi_*(G)$ is a torsion sheaf supported
firstly on those $b\in B$ corresponding to irreducible jumping
curves $C$, where it has length $h^1(G_{F_b})$, and secondly
on  those $b$ corresponding to reducible fibres $F_b=C_1\cup_p C_2$
such that $h^1(G_{F_b})\neq 0$ where its length is again
equal to this $h^1$. By property T, $E$ is almost balanced on all
reducible fibres $F_b$. This implies easily that
$F_b$ has at most a unique component, say $C_1$
of degree $d_1$, such that
$h^1(G_{C_1})\neq 0$, and in this case
$$h^0(G_{C_1})= 0$$
$$h^1(G_{C_2}(-p))= 0 $$  hence
$$h^1(G_{F_b})=h^1(G_{C_1}).$$
Now we have
$$E_{C_1}\simeq s\O(j)\oplus (r-s)\O(j-1)$$
with $j=\lceil\frac{ad_1}{r}\rceil$, and it is
immediate from this that $h^1(G_{C_1})\neq 0$
only if $t_I>j$, in which case 
$$h^1(G_{F_b})=h^1(G_{C_1})= rt_I-ad_1-r.$$
It follows that the total $h^1$ coming from reducible fibres
equals $h(t.,a.)$, so one side of GRR yields $-J_{d,E}(a.)-h(t.,a.).$\par
Now the other side of GRR generally equals
$$(r1_X+c_1(G)+\frac{1}{2}(c_1^2-2c_2)(G))(1_X-\frac{1}{2}K_X+\chi (\O_X)[pt])_2\tag 3.16$$
where $[pt]$ is a point and $_2$ denotes degree-2 part. Clearly
$$\chi (\O_X)=1-g. $$
 Next, the canonical class
$K_X$ was computed in [R2][R3] as
$$K_X= -2s_i+(2g-2-m_i)F+R_i$$
for any $ i$ (we take $i=1$), where $F$ is a fibre.
Given this, the computation of (3.16) is routine, 
yielding the formula (3.15).\qed

\enddemo\ls\ls
\subheading{Example 3.2} Take $n=2, d=4$ and let $E$ be the 
tangent bundle of $\P^2$, with Chern classes $a=3, b=3$. It is 
known that in this case $N_4=620, g=725.$ It is easy to compute 
that $$m_1=284, R_1.L=5220.$$ We take $D=7s_1$ and compute that in 
this case the $h^1$ contribution from the reducible curves is 
$$h(t.a.)=9180,$$ therefore finally $$J_{4,T_{\P^2}}=7944.$$ It is 
easy to see that on an irreducible rational quartic the tangent 
bundle cannot have $\O(4)$ as a direct summand so the jumping 
quartics all have splitting type $(7,5)$. It is also easy to see 
that a conic-pair $C_1\cup_pC_2$ with $C_1,C_2$ irreducible is 
jumping iff $C_1$ and $C_2$ are tangent at $p$. 
\ls\subheading{Example 3.3} Again for $n=2$, let $E$ be a 
semistable bundle with $c_1=0, c_2=m-1$ corresponding to $m$ 
general points in $\P^2$ (cf. [OSS]). Then we compute 
$$J_{d,E}=(m-1)N_d.$$ For $m=2, \Cal J_{d,E}$ is just a hyperplane 
section of $V_d$ corresponding to a certain point in $\P^2$. For 
$m=3,$ we get an interesting class of quadric sections of 
$V_d$.\ls \subheading{Remark 3.4} Note that Theorem 3.1 may be 
applied to the restriction of $E$ to a general linear subspace  
$A_0\subset\P^n.$ Hence the Theorem generalizes immediately to the 
case where one of the incidence conditions on the rational curve 
becomes containment in $A_0.$ \vfill\eject \ls\ls 
\centerline{\bf{REFERENCES}} \ls 

\item{[F]} W. Fulton: 'Intersection theory' Springer 1984.\ss
\item{[G]} A. Givental: 'Problems in quantum K-theory', talk
at Southern California Algebraic Geometry Seminar, 2/00.\ss
\item{[K]} J.Koll\'ar: 'Rational curves on algebraic varieties'
Springer 1996.\ss
\item{[OSS]} Ch. Okonek, M. Schneider, H. Spindler: 'Vector
bundles on projective spaces' Birkh\"auser 1980.\ss
\item{[P]} R. Pandharipande: 'The canonical class of $M_{0,n}(\P^r,d)$ and
enumerative geometry' IMRN 1997, no.4, 173-186.\ss
\item{[R1]} Z. Ran: 'Bend, break and count' Isr. J. Math 111 (1999),
109-124.\ss
\item{[R2]} Z. Ran: 'Bend, break and count II' 
(Math. Proc. Camb. Phil. Soc., to appear)\ss
\item{[R3]} Z. Ran: 'On the variety of rational space curves' Isr. J. Math
(to appear)\ss
\item{[S]}  E. Sernesi: ' On the existence of certain families of curves' Invent. math. 75
(1984), 25-57.
\ss
\enddocument